\documentstyle[12pt,leqno]{article}

\newfam\frakfam
\font\teneufm=eufm10     \textfont\frakfam=\teneufm
\font\seveneufm=eufm7    \scriptfont\frakfam=\seveneufm
\font\fiveeufm=eufm5     \scriptscriptfont\frakfam=\fiveeufm
\def\frak{\fam\frakfam \teneufm}

\newfam\msbfam
\font\tenmsb=msbm10      \textfont\msbfam=\tenmsb
\font\sevenmsb=msbm7     \scriptfont\msbfam=\sevenmsb
\font\fivemsb=msbm5      \scriptscriptfont\msbfam=\fivemsb
\def\Bbb{\fam\msbfam \tenmsb}

\font\titlebf=cmbx10  scaled \magstep2

\newcommand{\ie}{\mbox{\em i.e.}}
\newcommand{\fq}{\mbox{${\Bbb F}_q\,$}}
\newcommand{\fp}{\mbox{${\Bbb F}_p\,$}}
\newcommand{\clo}{\mbox{${\overline {\Bbb F}}_p\,$}}
\newcommand{\nq}{\mbox{${\Bbb F}_{q^n}\,$}}

\newcommand{\sn}{\mbox{${\frak S}_n\,$}}
\newcommand{\op}{\mbox{${\mathcal O}_P\,$}}

\newcommand{\pf}{\mbox{${\Bbb P}(F)\,$}}

\newcommand{\ith}{\mbox{$i^{th}\;$}}

\newcommand{\degdiff}{\mbox{\rm degDiff}}
\newcommand{\gal}{\mbox{\rm Gal}}

\newtheorem{theorem}{Theorem}[section]
\newtheorem{lemma}[theorem]{Lemma}
\newtheorem{proposition}[theorem]{Proposition}
\newtheorem{corollary}[theorem]{Corollary}
\newtheorem{definition}[theorem]{Definition}

\newtheorem{example}[theorem]{\sc Example}

\begin{document}
\font\titlebf=cmbx10 scaled \magstep2

\centerline{\titlebf Extensions of algebraic function fields with}
\centerline{\titlebf complete splitting of all rational places}
\centerline{{\bf Vinay Deolalikar} \footnote[1]{ This work was
part of the author's Ph.D. thesis and was supervised by the late
Prof. Dennis Estes.}}
\medskip
\medskip
\centerline{\sc Dedicated to the late Prof. Dennis Estes}
\centerline{Mathematics Subject Classification 14G05, 14G50}

\section{Introduction}
Let $F/K$ be an algebraic function field in one variable over a finite field of constants $K$, \ie, $F$ is a finite algebraic extension of $K(x)$ where $x \in F$ is transcendental over $K$. For simplicity, $K$ is assumed algebraically closed in $F$. Let $E$ be a finite separable extension of $F$. Let $N(E)$ and $g(E)$ denote the number of places of degree one, and the genus, respectively, of $E$. Let $[E:F]$ denote the degree of this extension.

The study of algebraic function fields with many places of degree one (also referred to as {\em rational} places) has received considerable attention in recent years. This was  fuelled initially by applications to coding theory that were first established by Goppa \cite{Gop1} in 1981, and later, as interesting mathematical objects in their own right, having connections to several well studied problems in arithmetic algebraic geometry. Many authors have also written on this subject in the language of curves over finite fields with many rational points.

In a previous paper \cite{DeoEst1}, we studied the problem of producing algebraic function fields with many places of degree one from a slightly different angle. We looked at these fields as extensions of the rational function field, and  described some families of such extensions in which all, or almost all (in most of our examples, that meant all except one) of the places of degree one in the rational function field split completely in the extension. We showed how many existing function fields that are known to be maximal in some sense (Hasse-Weil, Oesterle) actually fall in some of these families, and we presented more, hitherto unknown families with the same amount of splitting of rational places, but lower genus. In some cases, these families contained function fields that were maximal in the Oesterle sense. In that paper, we had used the notion of symmetry of  polynomials for an extension of finite fields $K/L$, to produce extensions of function fields over $K$ in which almost all rational places split completely.

In this paper, we introduce a larger class of functions, those that we call ``quasi-symmetric,'' that can be used to the same effect. The notion of quasi-symmetry generalises that of symmetry, for purposes of use in producing extensions in which almost all rational places split completely. We show how, using quasi-symmetric functions, it is actually possible to split {\em all} the rational places in extensions of function fields. Thus, some of the extensions described in this paper can attain the maximum possible value for the ratio of $N(E)/[E:F]$, for a fixed $[E:F]$.

Furthermore, these constructions are explicit, and therefore can be used for constructions of good codes. The class field theoretic descriptions \cite{Ser1,Ser2, Ser3} of  extensions in which many rational places split completely, do not provide explicit equations for the generators, and thus can not be used for applications to coding theory.

Another point that needs mention is that while most of the explicitly constructed extensions having lot of splitting of rational places have low genera and, therefore, few rational points, the methods described in this paper can be applied to construct extensions of arbitrarily high genera and number of rational points.

\section{Preliminaries}

Throughout this paper, we will use the following notation:

For symmetric polynomials:
\begin{tabbing}
\sn \hspace{1.6cm} \= the symmetric group on $n$ characters \\
${\bf s}_{n,i}(X)$ \> the \ith elementary symmetric polynomial on $n$ variables \\
$q$  \>  a power of a prime $p$ \\
${\Bbb F}_l$  \> the finite field of cardinality $l$ \\
$s_{n,i}(t)$ \> the \ith $(n,q)$-elementary symmetric polynomial
\end{tabbing}

For function fields and their symmetric subfields:
\begin{tabbing}
$K$ \hspace{1.6cm} \= the finite field of cardinality $q^n$, where $n>1$ \\
$F/K$ \> an algebraic function field in one variable whose full field of constants is $K$ \\
$F_s$ \> the subfield of $F$ comprising $(n,q)$-symmetric functions \\
$F_s^\phi$ \> the subfield of $F_s$ comprising functions whose coefficients are from \fq \\
$F_{qs}$ \> the subfield of $F$ comprising $(n,q)$-quasi-symmetric functions \\
$F^{\phi}_{qs}$ \> the subfield of $F_{qs}$ comprising functions whose coefficients are from \fq \\
$U_{qs}$ \> the \nq-vector space of $(n,q)$-quasi-symmetric functions from \nq to \nq\\
$V_{qs}$ \> the \fq-vector space of $(n,q)$-quasi-symmetric functions from \nq to \nq with a \\
    \> polynomial representation with \fq coefficients\\
$E$ \> a finite separable extension of $F$, $E=F(y)$ where $\varphi(y) = 0$
for some irreducible  \\
 \> polynomial $\varphi[T] \in F[T]$
\end{tabbing}

For a generic function field $F$:
\begin{tabbing}
${\Bbb P}(F)$  \hspace{1.3cm} \=  the set of places of $F$ \\
$N(F)$  \> the number of places of degree one in $F$ \\
$N_m(F)$\> the number of places of degree $m, m>1$, in $F$ \\
$g(K)$ \> the genus of $F$   \\
$P$ \> a generic place in $F$  \\
$v_P$ \> the normalized discrete valuation associated with the place $P$ \\
\op \> the valuation ring of the place $P$ \\
$P'$ \> a generic place lying above $P$ in a finite separable extension of $F$\\
$e(P'|P)$ \> the ramification index for $P'$ over $P$ \\
$d(P'|P)$ \> the different exponent for $P'$ over $P$
\end{tabbing}

For the rational function field $K(x)$:
\begin{tabbing}
$P_\alpha$ \hspace{1.6cm} \= the place in $K(x)$ that is the unique zero of $x-\alpha, \; \alpha \in K$ \\
$P_\infty$ \> the place in $K(x)$ that is the unique pole of $1/x$
\end{tabbing}

We also reproduce here some results from \cite{DeoThe} and \cite{DeoEst1}.

\begin{proposition} \label{proposition:Art-Sch}
Let $F/K$ be an algebraic function field, where $K=\nq$ is algebraically closed in $F$. Let $w \in F$ and assume that there exists a place $P \in \pf$ such that
$$ v_P(w) = -m, \, m > 0 \mbox{  and  } \gcd(m,q)  = 1.  $$
Then the polynomial $l(T)-w = a_{n-1}T^{q^{n-1}} + a_{n-2}T^{q^{n-2}}+ \ldots + a_0T - w \in F[T] $ is absolutely irreducible. Further, let $l(T)$ split into linear factors over $K$. Let $E=F(y)$ with
$$  a_{n-1}y^{q^{n-1}} + a_{n-2}y^{q^{n-2}}+ \ldots + a_0y = w. $$
Then the following hold:
\begin{enumerate}
\item $E/F$ is a Galois extension, with degree $[E:F] = q^{n-1}$. $\gal(E/F) = \{\sigma_\beta:y \rightarrow y + \beta\}_{l(\beta)=0}$.
\item $K$ is algebraically closed in $E$.
\item The place $P$ is totally ramified in $E$. Let the unique place of $E$ that lies above $P$ be $P'$. Then the different exponent $d(P'|P)$ in the extension $E/F$ is given by
$$ d(P'|P) = (q^{n-1}-1)(m+1).$$
\item Let $R \in \pf$, and $v_R(w) \geq 0$. Then $R$ is unramified in $E$.
\item If $a_{n-1}=\ldots=a_0=1$, and if $Q\in \pf$ is a zero of $w-\gamma$, with $\gamma \in \fq$. Then $Q$ splits completely in $E$.
\end{enumerate}
\end{proposition}
{\bf Proof}. For (i) - (iv), pl. refer \cite{Sul1}. For (v), notice that under the hypotheses, the equation $ T^{q^{n-1}} +T^{q^{n-2}}+ \ldots + T = \gamma $ has $q^{n-1}$ distinct roots in $K$.  \hfill $\Box$

For many of the extensions that we will describe, there exists no place where the hypothesis of Proposition~\ref{proposition:Art-Sch} is satisfied, namely, that the valuation of $w$ at the place is negative and coprime to the characteristic. In particular, we need a criterion for determining the irreducibility of the equations that we will need to use. We provide such a criterion in Proposition~\ref{proposition:irreducibility} and Corollary~\ref{corollary:irreducibility}.

\begin{proposition} \label{proposition:irreducibility}
Let $V$ be a finite subgroup of the additive group of \clo. Then $V$ is a \fp-vector space. Define $L_V(T) = \prod_{v \in V}(T-v)$. Thus, $L_V(T)$ is a separable \fp-linear polynomial whose degree is the cardinality of $V$. Now let $h(T,x) =  L_V(T) - f(x)$, where $f(x) \in \clo[x]$. Then, $h(T,x)=L_V(T) - f(x)$ is reducible over $\clo[T,x]$ iff there exists a polynomial $g(x) \in \clo[x]$ and a proper additive subgroup $W$ of $V$ such that $f(x) = L_{W'}(g(x))$, where $W' = L_W(V)$.
\end{proposition}
For a proof of this proposition, please refer to \cite{DeoThe} or \cite{DeoEst1}.

\begin{definition}
For $f(x) \in \clo[x]$, a coprime term of $f$ is a term with non-zero coefficient in $f$ whose degree is coprime to $p$. The coprime degree of $f$ is the degree of the coprime term of $f$ having the largest degree.
\end{definition}

\begin{corollary} \label{corollary:irreducibility}
Let $f(x) \in \clo[x]$. Let there be a coprime term in $f(x)$ of degree $d$, such that there are no terms of degree $dp^i$ for $i>0$ in $f(x)$. Then $L_V(T) -f(x)$ is irreducible for any subgroup $V \subset \clo$.
\end{corollary}
{\bf Proof}. Suppose $f(x)$ is the image of a linear polynomial $\sum a_nx^{p_n}$. Then the coprime term can only occur in the image of the term $a_0x$. But then, the images of the coprime term under $a_nx^{p^n}$, for $n>0$ will have degrees that contradict the hypothesis.

\begin{lemma} \label{lemma:subextensions}
 Let $F = K(x)$, where $K=\nq, q = p^m, r=m(n-1),$ and $E = F(y)$, where $y$ satisifes the following equation:
$$ y^{q^{n-1}} + y^{q^{n-2}} + \ldots + y = f(x), $$
and $f(x) \in F$ is not the image of any element in $F$ under a linear polynomial.
Then the following hold:
\begin{enumerate}
\item $E/F$ is a Galois extension of degree $[E:F]=q^{n-1}$. $\gal(E/F) = \{\sigma_\beta: y \rightarrow y + \beta \}_{s_{n,1}(\beta)=0}$ can be identified with the set of elements in \nq whose trace in \fq is zero by $\sigma_\beta \leftrightarrow \beta$. This gives it the structure of a $r$-dimensional \fp vector space.
\item There exists a tower of subextensions
$$ F=E^0 \subset E^1 \subset \ldots \subset E^{r} = E, $$
such that for $0 \leq i \leq r-1,\; [E^{i+1}: E^i]$ is a Galois extension of degree $p$.
\item Let $\{b_i\}_{1 \leq i \leq r}$ be a \fp-basis for $\gal(E/F)$. Then we can build the tower of subextensions as follows. We set $E^j$ to be  the fixed field of the subgroup of the Galois group that corresponds to the \fp-subspace generated by $\{b_1,b_2,\ldots, b_{r-j}\}$. Then, the generators of $E^{j}$ are $\{y_1,y_2,\ldots,y_j\}$, where $y_1,y_2,\ldots,y_{r}=y$ satisfy the following relations:
\begin{eqnarray*}
       y^p - B_{r}^{p-1} y &=& y_{r-1}, \\
       y_{r-1}^p - B_{r-1}^{p-1} y_{r-1} &=& y_{r-2}, \\
          \vdots &  & \vdots \\
       y_1^p - B_{1}^{p-1}y_1 &=& f(x),
\end{eqnarray*}
where,
$$
\begin{array}{rcll}
  \beta_{r,j} &=& b_{r-j+1},  & B_{r} = \beta_{r,r}, \\
  \beta_{r-1,j} &=& \beta_{r,j}^p - B_{r}^{p-1}\beta_{r,j}, & B_{r-1} = \beta_{r-1,r-1},\\
 \vdots & & \vdots &\vdots  \\
  \beta_{1,j}   &=& \beta_{2,j}^p - B_2^{p-1}\beta_{2,j}, &  B_1 = \beta_{1,1}.
\end{array}
$$
\end{enumerate}
\end{lemma}
For a proof of this lemma, please refer to \cite{DeoThe} or \cite{DeoEst1}.

\section{Quasi-symmetric functions}

We first recollect some notions about symmetric functions from \cite{DeoEst1}.

Let $R$ be an integral domain and $\overline{R}$ its field of fractions. Consider the polynomial ring  in $n$ variables over $R$, given by  $R\,[X]= R\,[x_1,x_2,\ldots,x_n]$. The symmetric group \sn acts in a natural way on this ring by permuting the variables.

\begin{definition}
A polynomial ${\bf f}(X) \in R\,[X]$ is said to be symmetric if it is fixed under the action of \sn. If \sn is allowed to act on $\overline{R}(X)$ in the natural way, its fixed points will be called symmetric rational functions, or simply, symmetric functions. These form a subfield $\overline{R}(X)_s$  of $\overline{R}(x)$. Furthermore, $\overline{R}(X)_s$ is generated by the $n$ elementary symmetric functions given by
\begin{eqnarray*}
  {\bf s}_{n,1}(X) &=& \sum_{i=1}^{n} x_i, \\
  {\bf s}_{n,2}(X) &=& \sum_{i<j \atop 1 \leq i,j \leq n} x_ix_j, \\
     \vdots  & &  \vdots \\
  {\bf s}_{n,n}(X) &=& x_1x_2\ldots x_n.
\end{eqnarray*}
\end{definition}

\begin{definition}
For the extension ${\Bbb F}_{q^n}/{\Bbb F}_q$, we will evaluate the elementary symmetric polynomials (resp. symmetric functions) in ${\Bbb F}_{q^n}(X)$ at $(X)=(t,\phi(t),\ldots,\phi^{n-1}(t))=(t,t^q,\ldots,t^{q^{n-1}})$. These will be called the $(n,q)$-elementary symmetric polynomials (resp. $(n,q)$-symmetric functions). For ${\bf f}(X) \in \nq\!(X)$, we will denote ${\bf f}(t,t^q,\ldots,t^{q^{n-1}})$ by $f(t)$, or, when the context is clear, by $f$.
\end{definition}

Thus the  $(n,q)$-elementary symmetric polynomials are the following:
\begin{eqnarray*}
s_{n,1}(t) &=& \sum_{0 \leq i \leq n-1} t^{q^i}, \\
s_{n,2}(t) &=& \sum_{i<j \atop 0 \leq i,j \leq n-1}t^{q^i}t^{q^j}, \\
    \vdots & & \vdots  \\
s_{n,n}(t) &=& t^{1+q+q^2+\ldots+q^{n-1}}.
\end{eqnarray*}

In \cite{DeoEst1}, we had demonstrated the use of $(n,q)$-symmetric functions in splitting places of degree one in extensions of algebraic functions fields.

In this paper, we extend the notion of symmetry to get a larger class of functions that can be very effectively used to split all places of degree one in extensions of function fields. We call these functions ``$(n,q)$-quasi-symmetric.''

\begin{definition} A polynomial ${\bf f}(X)$ in $R[X]$ will be called quasi-symmetric if it is fixed by the cycle $\varepsilon = (1\; 2\;\ldots\; n) \in {\frak S}_n$.  If $\varepsilon$ is allowed to act on $\overline{R}(X)$ in the natural way, its fixed points will be called quasi-symmetric rational functions, or simply, quasi-symmetric functions. These form a subfield $\overline{R}(X)_{qs}$  of $\overline{R}(X)$.
\end{definition}

\begin{lemma} For $n>2$, there always  exist quasi-symmetric functions that are not symmetric.
\end{lemma}
{\bf Proof}.  $\langle{\varepsilon}\rangle$ has index $(n-1)!$ in $\sn$. Thus for $n >2$, the set of functions fixed by \sn is strictly smaller than those fixed by $(\varepsilon)$. For $n=2,\; \sn = (\varepsilon)$ so that the notions of symmetric and quasi-symmetric coincide.\hfill $\Box$

\begin{example} \rm \label{example:quasi-symmetric-n=3}
\rm $(n=3)$ A family of quasi-symmetric functions in three variables is given below:
$${\bf f}(x_1,x_2,x_3) = x_1{x_2}^i + x_2{x_3}^i + x_3{x_1}^i. $$
Note that for $i \neq 0 \mbox{ or }1$, these are not symmetric.
\end{example}

\begin{definition}
Consider the  extension ${\Bbb F}_{q^n}/{\Bbb F}_q$ of finite fields. We will evaluate the  quasi-symmetric polynomials (resp. quasi-symmetric functions) in ${\Bbb F}_{q^n}(X)$ at  $(X)=(t,\phi(t),\ldots,\phi^{n-1}(t))=(t,t^q,\ldots,t^{q^{n-1}})$. These will be called $(n,q)$-quasi-symmetric  polynomials (resp. $(n,q)$-quasi-symmetric functions).
\end{definition}

\begin{example}\rm
Using the three-variable quasi-symmetric functions of Example~$\ref{example:quasi-symmetric-n=3}$, we can obtain  the following $(3,q)$-quasi-symmetric functions:
$$f(t) = {\bf f}(t,t^q,t^{q^2}) = t^{1+iq} + t^{q+iq^2} + t^{q^2 + i}.$$
Again, these are not symmetric for  $i \neq 0 \mbox{ or }1$.
\end{example}

\begin{lemma} \label{lemma:quasi-is-invariant-ifpart}
Let ${\bf f}(X) \in R[X]$. Then if ${\bf f}(X)$ is quasi-symmetric, $f(t^q) = f(t) \bmod (t^{q^n}-t)$.
\end{lemma}
{\bf Proof}.  We have that
$$ {\bf f}(\varepsilon(X)) = {\bf f}(x_2,x_3,\ldots,x_n,x_1). $$
Evaluating ${\bf f}(\varepsilon(X))$ at $(X)=(t,t^q,\ldots,t^{q^{n-1}})$, we get
$$  {\bf f}(t^q,t^{q^2},\ldots,t^{q^{n-1}},t) =  f(t^q)\bmod(t^{q^n} - t). $$
But since $f(X)$ is quasi-symmetric, this is equal to ${\bf f}(X)$ evaluated at $(X)=(t,t^q,\ldots,t^{q^{n-1}})$, which is ${\bf f}(t, t^q,\ldots,t^{q^{n-1}}) = f(t)$.   \hfill $\Box$

There is a converse to this statement if we allow only polynomials with degrees less than $q^n$. But first we need the following definition.

\begin{definition}
A lift of $f(t)$ is a polynomial ${\bf f}(X)$ such that ${\bf f}(t,t^q,\ldots,t^{q^{n-1}}) = f(t)$.
\end{definition}

\begin{lemma}  \label{lemma:quasi-is-invariant}
Let $f(t) \in \nq[t]$ have degree less than $q^n$. Then $f(t)$ is $(n,q)$-quasi-symmetric  iff $f(t^q) = f(t)\bmod (t^{q^n}-t)$.
\end{lemma}
{\bf Proof}. $\Rightarrow$ As in Lemma~\ref{lemma:quasi-is-invariant-ifpart}.

$\Leftarrow$ Each monomial in $f$ of degree $d$, with $d<q^n$ has a unique lift $x_1^{d_1}\ldots x_n^{d_n}$, where $d_i < q, 1 \leq i \leq n$ and $d_nd_{n-1}\ldots d_1$ is the base $q$ representation of $d$.

Then the lift of $f$ is just the sum of the lifts of each of its monomials. Now define $f_q(t)$ as the unique polynomial of degree less than $q^n$ which is congruent to $f(t^q)$ modulo $(t^{q^n}-t)$.

Now note that the lift of $f_q$ is just the cyclic shift in the variables of the lift of $f$. But by hypothesis, $f(t) = f_q(t)$ and thus the lift of $f$ must be invariant under a cyclic shift (\ie, under the action of $\varepsilon$), and thus is quasi-symmetric. It then follows that $f(t)$ is $(n,q)$-quasi-symmetric. \hfill $\Box$

\begin{corollary} \label{corollary:f-is-quasi}
A polynomial $f(t) \in \nq[t]$ is $(n,q)$-quasi-symmetric iff $f(\gamma^q) = f(\gamma), \;\; \forall \gamma \in \nq$.
\end{corollary}

\begin{lemma}
For a $(n,q)$-quasi-symmetric polynomial $f(t)$ with \fq coefficients, $f(\gamma) \in {\Bbb F}_q,\; \forall \gamma \in {\Bbb F}_{q^n}$.
\end{lemma}
{\bf Proof}. It suffices to prove that $\phi(f(\gamma))= f(\gamma)$. But this is immediate from the fact that the coefficients of $f$ are in ${\Bbb F}_q$ and  from Corollary~\ref{corollary:f-is-quasi}.  \hfill $\Box$

\begin{lemma}
Let $f$ be any function from \nq to \nq that satisfies $f(t^q) = f(t)$. Then  $f|_{\nq} = g(t)$, where $g(t)$ is a $(n,q)$-quasi-symmetric polynomial.
\end{lemma}
{\bf Proof}. By Lagrange interpolation, we can find a polynomial $g(t)$ in $\nq[t]$ that agrees with $f$ on \nq. Thus, $g(\gamma^q) = g(\gamma), \forall \gamma \in \nq$. Now from Corollary~\ref{corollary:f-is-quasi}, it follows that $g(t)$ is $(n,q)$-quasi-symmetric.\hfill $\Box$

Due to this result, we may use the terms $(n,q)$-quasi-symmetric function and $(n,q)$-quasi-symmetric polynomial interchangeably when the context allows.

\begin{definition}
Let $U_{qs}$  denote the ${\Bbb F}_{q^n}$-vector space of $(n,q)$-quasi-symmetric functions from \nq to $\nq\!,$ and $V_{qs} \subset U_{qs}$ denote the ${\Bbb F}_q$-space of all functions $f \in U_{qs}$ such that $f|_{\nq} =  g(t)|_{\nq}$ for some $g(t) \in {\Bbb F}_q[t]$. Thus, $V_{qs}$ consists of $(n,q)$-quasi-symmetric functions having a polynomial representation with coefficients in \fq.
\end{definition}

\begin{lemma}
Any ${\Bbb F}_q$-linearly independent subset $\{f_i(t)\}_{1 \leq i \leq r}$ in $V_{qs}$ is also a ${\Bbb F}_{q^n}$-linearly independent subset in $V_{qs}$.
\end{lemma}
{\bf Proof}. Let $u_1f_1(t) +\dots + u_rf_r(t) = 0, u_i \in {\Bbb F}_{q^n}$. Now express the $u_i$ as $\sum a_{ij}w_j, a_{ij} \in {\Bbb F}_q, w_1,\dots,w_n$ a ${\Bbb F}_q$-basis for ${\Bbb F}_{q^n}$, then $\sum(\sum a_{ij}f_i(t))w_j = 0$ and since $f_i(u) \in {\Bbb F}_q, u \in {\Bbb F}_{q^n}$, it follows that
$\sum a_{ij}f_i(t)=0$ and hence $a_{ij} = 0$.\hfill $\Box$

\begin{definition}
Let ${\mathcal O}$ denote the set of orbits of Galois for the action of $\gal({\Bbb F}_{q^n}/{\Bbb F}_q)$ on \nq.
\end{definition}

We note that $U_{qs}$ consists of functions from ${\Bbb F}_{q^n}$ to ${\Bbb F}_{q^n}$ which are constant on each orbit in ${\mathcal O}$. $V_{qs}$ consists of functions from \nq to \fq that are constant on these orbits.

\begin{lemma}
Let $f \in V_{qs}$, then there is a  unique polynomial representation $g(t) \in \fq[t]$, of degree $< q^n$, such that $f|_{\nq} = g(t)$.
\end{lemma}
{\bf Proof}. Choose a polynomial representation $g(t)$ of $f$ degree less that $q^n$. Define $\phi(g(t))_q$ as the unique polynomial of degree less than $q^n$ which is congruent to $\phi(g(t))$ modulo $(t^{q^n} - t)$. Since $f \in V_{qs}$,  we must have $g(t), \phi(g(t))_q$ agree on \nq. Therefore, they must be equal as their degree is less than $q^n$ (for if they were not, then the polynomial $g(t)-\phi(g(t))_q$ would have more zeros than its degree).  Thus, the coefficients of $g(t)$ are fixed by the Galois group and must lie in ${\Bbb F}_q$. Uniqueness follows similarly.\hfill $\Box$

When there is no scope for confusion, we will often say that a polynomial is in $V_{qs}$ to mean that it is the representative of a function in $V_{qs}$.

Now we examine the dimensions of $U_{qs}$ and $V_{qs}$.

\begin{theorem}
$U_{qs}$ and $V_{qs}$ have dimension $|{\mathcal O}|$ as vector spaces over \nq and \fq respectively
\end{theorem}
{\bf Proof}.  Consider the set of functions $\{f_i\}_{1 \leq i \leq |{\mathcal O}|}$, where $f_i$ is a function that is $1$ on the \ith orbit of Galois and $0$ on all others. Then this set is linearly independent over \nq (resp. \fq), and it spans $U_{qs}$ (resp. $V_{qs}$). \hfill $\Box$

\begin{lemma} \label{lemma:qs-no-zeros}
There exist $(n,q)$-quasi-symmetric functions that have no zeros in ${\Bbb
F}_{q^n}$.
\end{lemma}
{\bf Proof}. We are free to assign any value in \fq for a $(n,q)$-quasi-symmetric function on an orbit of Galois. In particular, we are free to assign a non-zero value to each orbit. Then there are $(q-1)^{|{\mathcal O}|}$ such functions constant on the orbits of Galois, and having no zero in \nq.\hfill $\Box$

We now discuss some ways to actually construct such functions whose existence is guaranteed by Lemma~\ref{lemma:qs-no-zeros}. We first consider a type of construction that we call the ``composition with irreducibles'' construction.

\begin{lemma} \label{lemma:composition-with-irreducibles}
Let $i(t)$ be a polynomial of degree $d_1$ over ${\Bbb F}_q$ that has no roots in \fq. In particular, we may choose $i(t)$ to be irreducible over \fq. Let $s(t) \in V_{qs}$ be a $(n,q)$-quasi-symmetric polynomial of degree $d_2$. Then the following hold:
\begin{enumerate}
\item $i(s(t))$ is $(n,q)$-quasi-symmetric, for $1 \leq i \leq n$, and maps ${\Bbb F}_{q^n}$ to ${\Bbb F}_q$.
\item $i(s(t))$ has no zeros in ${\Bbb F}_{q^n}$.
\end{enumerate}
\end{lemma}
{\bf Proof}.
For (i), note that $s(t)$  maps ${\Bbb F}_{q^n}$ to ${\Bbb F}_q$. Now since the coefficients of $i(t)$ are from \fq, we get the result.  For (ii), assume the contrary, \ie, let
there be $\alpha \in {\Bbb F}_{q^n}$ s.t $f(s(\alpha)) = 0$. Then we have that $s(\alpha)$ is a zero of $i(t)$. But $s(\alpha) \in {\Bbb
F}_q$, since $s(t) \in V_{qs}$. This would then give us a zero of $i(t)$ in ${\Bbb F}_q$, contradicting the hypothesis. \hfill $\Box$

The resulting polynomial obtained by this composition would have degree equal to $d_1d_2$. Since there exist irreducibles of any arbitrary degree $>1$ over \fq, we can now construct $(n,q)$-quasi-symmetric polynomials that map \nq to \fq  having degree equal to any multiple of a $(n,q)$-quasi-symmetric function in $V_{qs}$, and without a zero in ${\Bbb F}_{q^n}$.

A particularly useful form of this construction is given below.

\begin{corollary}
Let $s(t) \in V_{qs}$ be a $(n,q)$-quasi-symmetric polynomial. Then the polynomial $f(t)= {s(t)}^m - \beta$, where $\beta \in {\Bbb F}_q$ is not a $m^{th}$ power in ${\Bbb F}_q$, is $(n,q)$-quasi-symmetric,
maps ${\Bbb F}_{q^n}$ to ${\Bbb F}_q$, and has no zeros in ${\Bbb F}_{q^n}$.
\end{corollary}
{\bf Proof}. Follows from Lemma~\ref{lemma:composition-with-irreducibles}, by choosing $i(t) = t^m - \beta$.\hfill $\Box$

\begin{example}\rm
The polynomial $f(t) = t^{10} + 2t^6 + t^2 - 2$ is $(2,5)$-quasi-symmetric and has no zeros in ${\Bbb F}_{25}$. We obtain $f(t)$ by composing the irreducible $ i(t) = t^2 - 2$ with the first $(2,5)$-elementary symmetric polynomial, given by $s_{2,1}(t) = t^5 + t$.
\end{example}

Now we wish to use these notions in the setting of the algebraic function field $F/K$, where $K =\nq$.
\begin{definition}
$F_{qs}$ will denote  the field of $(n,q)$-quasi-symmetric functions in $F$. $F^{\phi}_{qs}$ will denote the field of $(n,q)$-quasi-symmetric rational functions in $F$, whose coefficients are from \fq.
\end{definition}

\section{Quasi-symmetric extensions of function fields}

$F$ and $K$ are as described earlier. $E$ is a finite separable extension of $F$, generated by $y$, where $\varphi(y) = 0$, for $\varphi(T)$ an irreducible polynomial in $F[T]$.

In this section we will introduce families of extensions of $F$ whose generators satisfy explicit equations involving only $(n,q)$-quasi-symmetric functions. Let $y$ satisfy
                        $$ g(y) = f(x),$$
where $f, g \in F_{qs}^\phi$. If $K=\nq$, this implies that in the residue field of a rational place, alhough the class of $x$ and $y$ can assume any values in ${\Bbb F}_{q^n} \cup \infty$, that of $f(x)$ and $g(y)$ will assume values only in ${\Bbb F}_{q}\cup \infty$. Among the Galois extensions that such equations can produce are the two special cases of extensions of Artin-Schreier and Kummer type.

\subsection{Quasi-symmetric extensions of Artin-Schreier type}

We are now in a position to state the main theorems of this section.
\begin{theorem} \label{theorem:general-quasi-symmetric}
Let $F=K(x)$ where $K=\nq$. Let $E=F(y)$, where $y$ satisfies the equation
\begin{equation} \label{equation:general-quasi-symmetric}
      y^{q^{n-1}} + y^{q^{n-2}} + \ldots + y = \frac{h(x)}{g(x)},
\end{equation}
where $h(x),g(x) \in F_{qs}^\phi$, and $\frac{h(x)}{g(x)}$ is not the image of any rational function in $F$ under a linear polynomial. Then the following hold:
\begin{enumerate}
\item $E/F$ is a Galois extension, with degree $[E:F] = q^{n-1}$. $\gal(E/F)= \{\sigma_\beta:y \rightarrow y + \beta\}_{s_{n,1}(\beta) = 0}.$
\item Let $P \in \pf$ be such that $v_P \left(\frac{h(x)}{g(x)}\right) = -m, \;m>0$ and $\gcd(m,q) = 1$. Then $P$ is totally ramified in $E$, with different exponent
$$ d(P'|P) = (q^{n-1} -1)(m + 1). $$
\item Let $Q \in \pf$ be any rational place, such that  $v_Q \left(\frac{h(x)}{g(x)}\right) \geq 0$. Then $Q$ splits completely in $E$.
\end{enumerate}
\end{theorem}
{\bf Proof}. For (i), irreducibility follows from Proposition~\ref{proposition:irreducibility}. The proof of that proposition carries through even for the case of rational functions. Also, we note that since  $h(x),g(x) \in F_{qs}^\phi$, the residue class of $\frac{h(x)}{g(x)}$ at any rational place is in \fq. Now (iii) follows from  Proposition~\ref{proposition:Art-Sch}.\hfill $\Box$

\begin{example}\rm
Let $F=K(x)$ and $K = {\Bbb F}_{q^3}$. Let $E=F(y)$, where $y$ satisfies the equation
$$ y^{q^2} + y^q + y = x^{1+iq} + x^{q +iq^2 } + x^{q^2+i}. $$
Now from Lemma~$\ref{lemma:subextensions}$, we can get a subextension $E^1$ whose degree over $F$ is $q$. Thus, $E^1 = F(y_1)$ where $y_1$ satisfies the equation
$$ y_1^q - b^{1-q}y_1 = x^{1+iq} + x^{q +iq^2 } + x^{q^2+i}, $$
where $b$ is a nonzero element of ${\Bbb F}_{q^3}$ whose trace in \fq is zero.
By substituting $y_1 = z_1 + x^{1+iq},$ we get the following equation
$$ z_1^q - b^{1-q}z_1 = (1+b^{1-q})x^{1+iq} + x^{q^2+i}. $$
In the case where $0 < i < q$, the degree of the RHS is $q^2 + i$. Also let $\gcd(i,q)=1$. Then we get that if $P_\infty^1$ is the place dividing $P_\infty$ in $E^1$,
$$ d(P_\infty^1|P_\infty) = (q -1)(q^2+i+1). $$
For the case of $i > q+1$, the degree of the RHS is $1+iq$  and then
$$ d(P_\infty^1|P_\infty) = (q -1)(2 + iq). $$
In case $i=q+1$, both the terms in the RHS have equal degrees. Notice that $1 + b^{1-q} = -b^{q^2-q} \neq 1$ unless the characteristic is three. Thus, if the characteristic is not three, we have that
$$ z_1^q - b^{1-q}z_1 = (1- b^{q^2-q})\;x^{q^2+q+1} \neq 0.$$
The different exponent is given by
$$ d(P_\infty^1|P_\infty) = (q^2 -1)(q^2+q+2). $$
In all these cases, we have that $N(E) = q^5+1$, since all rational places, except for $P_\infty$ split completely. For $i=1$ the RHS is the  $(3,q)$-elementary symmetric polynomial $s_{3,2}$. In that case, for $q=2$, $g(E) = 6$ and the extension attains the Oesterle lower bound on genus.
\end{example}

We now use $(n,q)$-quasi-symmetric functions that have no zeros in $K$, to build extensions of $F/K$.

\begin{theorem} \label{theorem:all-places-split}
Let $F=K(x)$, where $K = \nq$. Let $E=F(y)$, where $y$ satisfies the equation
 $$y^{q^{n-1}} + y^{q^{n-2}} + \ldots + y = \frac{h(x)}{g(x)}.$$
where $h(x), g(x) \in F_{qs}^\phi$, $\frac{h(x)}{g(x)}$ is not the image of a rational function in $F$ under a linear polynomial, $deg(g(x)) > deg(h(x))$, and $g(x)$ has no zeros in \nq. Then all the rational places of $F$ split completely in $E$. Thus $N(E) = q^{n-1}(q^n +1)$, and $N(E)/[E:F]$ attains its maximum possible value of $q^n + 1$.
\end{theorem}
{\bf Proof}. For all $\alpha \in \nq$, $ h(\alpha)/g(\alpha) \in \fq$, which ensures splitting of all rational places of the form $P_\alpha$. Also, since $deg(g(x)) > deg(h(x))$, the RHS of the equation has a zero at $P_\infty$, ensuring that $P_\infty$ also splits completely in $E$.\hfill $\Box$

\begin{theorem}
Consider the extension $E/F$ of Theorem$~\ref{theorem:all-places-split}$. We can find subextensions
$$ F=E^0 \subset E^1 \subset \ldots \subset E^{r} = E, $$
such that for $0 \leq i \leq r-1,\; [E^{i+1}: E^i]$ is a Galois extension of degree $p$, and in each extension $E^i/F$, all the rational places split completely.
\end{theorem}
{\bf Proof}. Refer to Lemma~\ref{lemma:subextensions}.

\begin{example} \rm
Let $F=K(x)$, where $K = {\Bbb F}_{q^2},\; p \neq 2$. Let $E=F(y)$, where $y$ satisfies the equation
$$y^q + y = \frac{x^{q+1}}{x^{2q} + 2x^{q+1} + x^2 - \alpha},$$
where $\alpha$ is not a square in \fq. All the rational places in $F$ split completely in $E$. The non-rational places which divide $x^{2q} + 2x^{q+1} + x^2 - \alpha$ are totally ramified.
For the special case of $q=5$, and obtaining a quadratic extension of ${\Bbb F}_5$ using the root $t$ of the irreducible polynomial $y^2 - 3$, we get the following  factorization into irreducibles of $x^{10} + 2x^{6} + x^2 - 2 $, where $2$ is our chosen non-square in  ${\Bbb F}_5$:
$$ x^{10} + 2x^{6} + x^2 - 2 = (x^5 + x + 2t)(x^5 + x + 3t). $$
This then gives us the following degree of different for $E/F$
$$ \degdiff(E/F) = 2(5-1)(1+1)(5) = 80, $$
which gives us the genus $g(E) = 37$. Also, since all rational places split completely, $N(E) = 130$. For this value of $q$ and $N$, the Oesterle lower bound on genus is $g \geq 11$. Thus while all the rational places split completely in this extension, the rise in genus that we incur in achieving this maximum splitting is high. Thus the $N/g$ ratio for $E$ is not very high. This is a typical example in this respect. Note that we can split all except one rational place of ${\Bbb F}_{25}$ in a degree $5$ extension by using the Hermitian function field, given by $E = F(y)$, where $y$ satisfies the equation
$$ y^5 + y = x^{6}.$$
In this case, the only place in \pf that is ramified is $P_\infty$, which is totally ramified with different exponent given by
$$ d(P'_\infty|P_\infty) = 5+2 = 7, \mbox{  and,} $$
$$ \degdiff(E/F) = (5-1)(6+1) = 28.$$
Thus, to split the one remaining rational place, we have had to almost triple the degree of the different.
\end{example}

\subsection{Quasi-symmetric extensions of Kummer type}

We now study extensions whose Galois group is a subgroup of the
multiplicative group $K^*$. For this we will need that the field
contain a primitive \ith root of unity $\xi_i$ for some $i$ coprime to $p$. In particular we know that $K$ contains $\xi_i$ for $i = \frac{q^n-1}{q-1}$.

\begin{theorem} \label{theorem:Kummer-quasi-symmetric}
Let $F=K(x)$ where $K = \nq$. Let $E=F(y)$, where $y$ satisfies the equation
               $$ y^{\frac{q^n-1}{q-1}} = \frac{h(x)}{g(x)}, $$
where  $h(x), g(x) \in  F_{qs}^\phi$ and $\frac{h(x)}{g(x)} \neq w^{\frac{q^n-1}{q-1}}, \; \forall w \in F$. Then the following hold:
\begin{enumerate}
\item This is a cyclic Galois extension, with degree $[E:F] = \frac{q^n-1}{q-1}$. $\gal(E/F) = \{ \sigma_j : y \rightarrow y{\xi}^j \}_{1 \leq j \leq \frac{q^n-1}{q-1}}$.
\item The only places of $F$ that may be ramified are the unique pole $P_{\infty}$ of $x$ and the zeros of $h(x)$ and $g(x)$. For any such place $P$, let $v_P = v_{P}(h(x)/g(x))$ be the corresponding valuation. Let $P'$ denote a generic place lying above $P$ in $E$. Then the ramification index for $P'$ over $P$ is given by
$$ e(P'|P) =\frac{[E:F]}{r_{P}},  $$
where $r_{P} = \gcd([E:F],v_P) > 0$. Since the extension is tame, the different exponent for $P'$ over $P$ is given by
$$  d(P'|P) = e(P'|P)-1. $$
\item All other rational places of ${\Bbb F}_{q^n}(x)$ split completely in the
extension. The rational places discussed in {\rm (ii)}, if not ramified, also split completely.
\end{enumerate}
\end{theorem}
{\bf Proof}. The proofs of (i) and (ii) can be looked up in any
reference on Kummer extensions. An excellent reference is \cite{Sti1}. (iii) follows from the fact that for $\alpha \in {\Bbb F}_{q^n}$ and  not  a zero of $h(x)$ or $g(x)$, $h(\alpha)/g(\alpha) \in {\Bbb F}_q^*$ and then we know that it has  $\frac{q^n-1}{q-1}$ pre-images under the norm map.\hfill $\Box$

\begin{corollary}
If in Theorem~$\ref{theorem:Kummer-quasi-symmetric}$, $h(x)$ and $g(x)$ have no zeros in ${\Bbb F}_{q^n}$, then other than possibly $P_{\infty}$, all rational places split completely in the extension, giving atleast $\frac{q^{2n} -q^n}{q-1}$ rational places. In addition, if $r_{\infty} = [E:F]$, then $P_{\infty}$ also splits completely. In that case, all rational places split completely in the extension.
\end{corollary}

\begin{example} \rm Let $q=p^e, e>1, p \neq 2, m= \frac{q-1}{p-1}$ and $\alpha$ not a square in ${\Bbb F}_p$. Let $F=K(x)$, where $K = \fq$, and $E=F(y)$, where $y$ satisfes the equation
 $$ y^m = x^{2m} - \alpha. $$
Now all the rational places will spilt completely. Ramification will be restricted to the non-rational places that correspond to the irreducible polynomials
which divide  $x^{2m} - \alpha$. Thus, this extension will attain the maximum possible value for the ratio $N(E)/[E:F]$, for $[E:F]=m$.
\end{example}

We have discussed a few techniques that can be used to split ``almost all'' rational places in extensions of function fields. In most of our examples, that means all, or all except one. We observe that it is very hard to keep the genus low if we split all rational places. It seems that to split the last rational place entails a fairly high increase in genus, as compared to splitting all but one rational places.

\section*{Acknowledgements}
I would like to express my deep sense of gratitude to Prof. Dennis Estes, who supervised this work,  and tragically passed away just prior to its completion. Without his constant help, I could not have made any progress whatsoever. This work is dedicated to him.

I would also like to profusely thank Joe Wetherell for all his help. He helped me immensely in completing this work following the demise of Prof. Dennis Estes.

\end{document}